%% file: L-series4.tex
 \documentclass[10pt,twoside,english]{article}
\usepackage{amssymb,amsmath,babel,geometry,latexsym,graphics,tabularx,shapepar,enumerate}
\input{option_keys}
\usepackage{hyperref}
\usepackage[all,2cell]{xy} \UseAllTwocells \SilentMatrices
\def\carlitz{{{\text{Car}}}}

\newtheorem{Theorem}{Theorem}
\newtheorem{Lemme}[Theorem]{Lemma}
\newtheorem{Proposition}[Theorem]{Proposition}
\newtheorem{Corollaire}[Theorem]{Corollary}
\newtheorem{Definition}[Theorem]{Definition}

\newcommand{\ol}[1]{\chi_t(#1)}

\newcommand{\ZZ}{\mathbb{Z}}
\newcommand{\FF}{\mathbb{F}}
\newcommand{\CC}{\mathbb{C}}
\newcommand{\QQ}{\mathbb{Q}}

\newcommand{\LL}{\mathbb{L}}
\newcommand{\MM}{\mathbb{M}}

\newcommand{\TT}{\mathbb{T}}

\newcommand{\bfs}{\boldsymbol{s}}

\newcommand{\bsb}{\boldsymbol}

\newcommand{\sqm}[4]{
\left(\begin{array}{ll}#1 & #2 \\ #3 & #4\end{array}\right)}
 
\newcommand\CVD{{\hfill\hfil{\lower 2 pt\hbox{\vrule\vbox to 7pt 
{\hrule width 6pt\vfill\hrule}\vrule}}}\vskip 0.5cm}

\let\ge=\geqslant                                 

\title{Values of certain $L$-series \\
in positive characteristic}
\author{Federico Pellarin\footnote{Current address: LaMUSE, 23, rue du Dr. Paul Michelon, 42023 Saint-Etienne Cedex.} \footnote{Supported by the contract ANR ``HAMOT", BLAN-0115-01.}}

\begin{document}

\maketitle

\begin{small}
\noindent\textbf{Abstract.} We introduce a family of $L$-series specialising to both
$L$-series associated to certain Dirichlet characters over $\FF_q[\theta]$ and to integral values of Carlitz-Goss zeta function associated to $\FF_q[\theta]$. We prove, with the use of the theory of deformations of vectorial modular forms, a formula for their value at $1$, as well as some arithmetic properties of values at other positive integers. 
\end{small}


%

\medskip


\section{Introduction, results}

\medskip

let $q=p^e$ be a power of a
prime number $p$ with $e>0$ an integer, let $\FF_q$ be the finite field with $q$ elements.
We consider, for an indeterminate $\theta$, the polynomial ring $A=\FF_q[\theta]$ and its fraction field $K=\FF_q(\theta)$. On $K$, we consider the absolute value
$|\cdot|$ defined by $|a|=q^{\deg_\theta a}$, $a$ being in $K$, so that
$|\theta| = q$.  Let $K_\infty :=\FF_q((1/\theta))$ be the
completion of $K$ for this absolute value, let $K_\infty^{\text{alg}}$ be an
algebraic closure of $K_\infty$, let $\CC_\infty$ be the completion of
$K_\infty^{\text{alg}}$ for the unique extension of $|\cdot|$ to $K_\infty^{\text{alg}}$, and let $K^{\text{sep}}$ be the separable closure of $K$
in $K_\infty^{\text{alg}}$.

We consider an element $t$ of $\CC_\infty$. We have the ``evaluating at $t$" ring homomorphism
$$\chi_t:A\rightarrow \FF_q[t]$$
defined by $\chi_t(a)=a(t)$. In other words, $\chi_t(a)$ is the image of the polynomial map $a(t)$ obtained by substituting, in $a(\theta)$, $\theta$ by $t$. For example,
$\chi_t(1)=1$ and $\chi_t(\theta)=t$. If we choose $t\in\FF_q^{\text{alg}}$ then $\chi_t$ factors through a Dirichlet character modulo the ideal generated by the minimal polynomial of $t$
in $A$. 

We can also consider $t$ as an indeterminate; 
for $\alpha>0$ an integer, we then have a formal series
$$L(\chi_t,\alpha)=\sum_{a\in A^+}\chi_t(a)a^{-\alpha}=\prod_{\mathfrak{p}}(1-\chi_t(\mathfrak{p})\mathfrak{p}^{-\alpha})^{-1},$$ 
where $A^+$ denotes the set of monic polynomials of $A$ and where the eulerian product runs over the monic irreducible 
polynomials of $A$, which turns out to be well defined in $K_\infty[[t]]$. 
This formal series 
converges for $\log_q|t|<\alpha$, $\log_q$ being the logarithm in base $q$. In this paper, we are
interested in the function $L(\chi_t,\alpha)$ of the variable $t\in\CC_\infty$, for a fixed positive integer $\alpha$.

We give some relevant examples of values of these series.
If $t=\theta$ and $\alpha>1$, then $L(\chi_\theta,\alpha)$ converges to the value of Carlitz-Goss' zeta series at $\alpha-1$:
$$L(\chi_\theta,\alpha)=\zeta(\alpha-1)=\sum_{a\in A^+}a^{1-\alpha}.$$
If on the other side we consider $t\in\FF_q^{\text{alg}}$, then, for $\alpha>0$, 
$L(\chi_t,\alpha)$ converges to the value at $\alpha$ of the $L$-series
associated to a Dirichlet character, see Goss' book \cite{Go2}.

We will need some classical functions related to Carlitz's module. First of all, its {\em exponential function} $e_{\carlitz}$ defined, for all $\eta\in \CC_\infty$, by the sum of the convergent series:
\begin{equation}\label{exponential}
e_\carlitz(\eta)=\sum_{n\ge 0}\frac{\eta^{q^n}}{d_n},
\end{equation}
where $d_0:=1$ and $d_i:=[i][i-1]^q\cdots[1]^{q^{i-1}}$, with $[i]=\theta^{q^i}-\theta$ if $i>0$.

We choose once and for all a fundamental period $\widetilde{\pi}$ of $e_{\carlitz}$.
It is possible to show that $\widetilde{\pi}$ is equal, up to the choice of a $(q-1)$-th root of $-\theta$,
to the (value of the) convergent infinite product:
$$\widetilde{\pi}:=\theta(-\theta)^{\frac{1}{q-1}}\prod_{i=1}^\infty(1-\theta^{1-q^i})^{-1}\in (-\theta)^{\frac{1}{q-1}}K_\infty.$$
Next, we need the following series of $K^{\text{sep}}[[t]]$:
\begin{equation}\label{scarlitz}
s_{\carlitz}(t):=\sum_{i=0}^\infty e_\carlitz\left(\frac{\widetilde{\pi}}{\theta^{i+1}}\right)t^i=\sum_{n=0}^\infty\frac{\widetilde{\pi}^{q^n}}{d_n(\theta^{q^n}-t)},
\end{equation} converging for $|t|<q$. 
We shall prove:
\begin{Theorem}\label{corollairezeta11}
The following identity holds:
$$L(\chi_t,1)=-\frac{\widetilde{\pi}}{(t-\theta) s_{\carlitz}}.$$ 
\end{Theorem}
The function $s_\carlitz$ is the {\em canonical rigid analytic trivialisation} of the so-called {\em Carlitz's motive}, property leading to three important facts that we recall now. First of all, $s_\carlitz$ generates the one-dimensional $\FF_q(t)$-vector space of solutions of the $\tau$-difference equation
\begin{equation}\label{scarlitz}
\tau X =(t-\theta)X,
\end{equation} in the fraction field of the ring of formal power series in $t$ converging for 
$|t|\leq 1$,
where $\tau:\CC_\infty((t))\rightarrow\CC_\infty((t))$ is the operator defined by $\tau\sum c_it^i=\sum c_i^qt^i$.
Second, $\bsb{\Omega}:=(\tau s_\carlitz)^{-1}$ is entire holomorphic on $\CC_\infty$ with the only zeros at
$t=\theta^q,\theta^{q^2},\ldots$ and 
third, we have the limit $\lim_{t\rightarrow\theta}\bsb{\Omega}(t)=-\frac{1}{\widetilde{\pi}}$.
We refer to \cite{Bourbaki} for a description of the main properties of
the functions $s_\carlitz,\bsb{\Omega}$, or to the papers \cite{An, ABP,Pa}, where $\bsb{\Omega}$ was originally introduced and studied.

The function $\bsb{\Omega}$ being entire, the function $L(\chi_t,1)=-\widetilde{\pi}\bsb{\Omega}(t)$ allows, beyond its domain of convergence, entire analytic continuation in terms of the variable $t$ with the only zeros at the points $t=\theta^q,\theta^{q^2},\ldots$ Also, it is well known (see \cite{ABP})
that $\bsb{\Omega}$ has the following infinite product expansion:
$$\bsb{\Omega}(t)=-\theta^{-1}(-\theta)^{-\frac{1}{q-1}}\prod_{i=1}^\infty(1-t\theta^{-q^i})\in(-\theta)^{\frac{1}{q-1}}K_\infty[[t]],$$ so that the product $\widetilde{\pi}\bsb{\Omega}$ no longer depends on a choice of $(-\theta)^{1/(q-1)}$. We deduce that the 
function $L(\chi_t,1)$, apart from its eulerian product, also has the following product expansion converging this time for all $t\in\CC_\infty$:
\begin{equation}\label{infiniteprod}
L(\chi_t,1)=\prod_{i=1}^\infty\frac{1-\frac{t}{\theta^{q^i}}}{1-\frac{\theta}{\theta^{q^i}}}=\prod_{i=1}^\infty\left(1-\frac{t-\theta}{[i]}\right).\end{equation}
It follows, from the above product expansion, that
\begin{equation}\label{limit1}
\lim_{t\to\theta}L(\chi_t,1)=1.\end{equation} This value coincides 
with $\zeta(0)$. Moreover, we obtain from (\ref{infiniteprod}) that $\lim_{t\to\theta^{q^k}}L(\chi_t,1)=0$ for 
$k>0$. This corresponds to the value at $1-q^k$ of $\zeta$. Indeed, we know more generally that $\zeta(s(q-1))=0$ for negative values of $s$. We see here a sign of the existence of functional equations for $\zeta$, see Goss' paper \cite{Go25} on ``zeta phenomenology". Indeed, in Theorem \ref{corollairezeta11}, the function $\tau s_\carlitz$ plays a role analogue to that of Euler's gamma function in the classical functional equation of Riemann's zeta function; its poles provide the trivial zeroes of  
Carlitz-Goss' $\zeta$ at the points $1-q^k$, $k>0$ and this is the first known interplay between certain
values of $\zeta$ at positive and negative multiples of $q-1$ (\footnote{This function was used in an essential way 
in \cite{ABP} in the study of the algebraic relations between values of the ``geometric" gamma function
of Thakur.}).

The above observations also agree with the results of Goss in \cite{Go3}, where he considers, for $\beta$ an integer, 
an analytic continuation of the function $\alpha\mapsto L(\chi_t^\beta,\alpha)$
with 
$$L(\chi_t^\beta,\alpha)=\sum_{a\in A^+}\chi_t(a)^\beta a^{-\alpha}$$ to the space $\mathbb{S}_\infty=\CC_\infty^\times\times \ZZ_p$ via the map $\alpha\mapsto((1/\theta)^{-\alpha},\alpha)\in\mathbb{S}_\infty$.

It is interesting to notice that Theorem \ref{corollairezeta11} directly implies the classical formulas for the values $$\zeta(q^k-1)=\sum_{a\in A^+}a^{1-q^k}=(-1)^k\frac{\widetilde{\pi}^{q^k-1}}{[k][k-1]\cdots[1]}$$ of Carlitz-Goss' zeta value at $q^k-1$ for $k>0$. This follows easily from the computation of the limit $t\rightarrow\theta$ in the
formula $\tau^kL(\chi_t,1)=-\widetilde{\pi}^{q^k}/(\tau^k(t-\theta)s_{\carlitz})$, observing that 
\begin{equation}\label{formula11}
\tau^k((t-\theta) s_\carlitz)=(t-\theta^{q^k})\cdots(t-\theta^q)(t-\theta) s_\carlitz\end{equation} by (\ref{scarlitz}).
 
Also, 
if $t=\xi\in\FF_q^\text{alg}$, Theorem \ref{corollairezeta11} implies that $L(\chi_\xi,1)$, the value of an $L$-function associated to a Dirichlet character,
is a multiple of $\widetilde{\pi}$ by an algebraic element of $\CC_\infty$. More precisely, we get the 
following result.

\begin{Corollaire}\label{fq}
For $\xi\in\FF_{q}^{\text{alg}}$ with $\xi^{q^r}=\xi$ ($r>0$), we have
$$L(\chi_\xi,1)=-\frac{\widetilde{\pi}}{\rho_\xi},$$ where 
$\rho_\xi\neq0$ is the root $(\tau s_\carlitz)(\xi)\in\CC_\infty$ of the polynomial equation
$$X^{q^r-1}=(\xi-\theta^{q^{r}})\cdots(\xi-\theta^q).$$
\end{Corollaire}
The corollary above follows from Theorem \ref{corollairezeta11} with $t=\xi$, applying (\ref{formula11}) with $k=r$,
simply noticing that $(\tau^r(\tau s_\carlitz))(\xi)=((\tau s_\carlitz)(\xi))^{q^r}$.

Some of the above consequences are also covered by the so-called Anderson's {\em log-algebraic power series} identities for {\em twisted harmonic sums}, see \cite{Anbis,Anter,Lutes}, see also \cite{Dam}. One of the new features
of the present work is to highlight the fact that several such results on values at one of $L$-series in positive
characteristic belong to the special family described in Theorem \ref{corollairezeta11}.

We briefly mention on the way that the series $\zeta(\alpha),L(\chi_t,\alpha)$ can be viewed as a special case of the (often divergent) multivariable formal series:
$$\mathcal{L}(\chi_{t_1},\ldots,\chi_{t_r};\alpha_1,\ldots,\alpha_r)=\sum_{a\in A^+}\chi_{t_1}^{-\alpha_1}(a)\cdots\chi_{t_r}^{-\alpha_r}(a).$$ For $r=1$ we get $\mathcal{L}(\chi_\theta,\alpha)=\zeta(\alpha)$ and 
for $r=2$ we get $\mathcal{L}(\chi_t,\chi_\theta;-1,\alpha)=L(\chi_t,\alpha)$. These series may also be of some interest
in the arithmetic theory of function fields of positive characteristic.

\medskip

For other positive values of $\alpha\equiv 1\pmod{q-1}$, we have a result on $L(\chi_t,\alpha)$
at once more general and less precise.

\begin{Theorem}\label{generalalpha}
Let $\alpha$ be a positive integer such that $\alpha\equiv1\pmod{q-1}$. There exists a non-zero element 
$\lambda_\alpha\in\FF_q(t,\theta)$ such that 
$$L(\chi_t,\alpha)=\lambda_\alpha\frac{\widetilde{\pi}^\alpha}{(t-\theta)s_\carlitz}.$$
\end{Theorem}
Theorem \ref{corollairezeta11} implies that $\lambda_1=-1$. On the other hand, again by (\ref{formula11}), we have the formula 
$$\tau^k\lambda_\alpha=(t-\theta^{q^k})\cdots(t-\theta^q)\lambda_{q^k\alpha}.$$
Apart from this, the explicit computation of the $\lambda_\alpha$'s is difficult and very little is known
on these coefficients which could encode, we hope, an interesting generalisation in $\FF_q(t,\theta)$ of the theory of Bernoulli-Carlitz numbers.

\medskip

\noindent\emph{Transcendence properties.} For $t$ varying in $K^{\text{alg}}$, algebraic closure of $K$ in $\CC_\infty$, our results can be used to obtain transcendence properties
of values $L(\chi_{t},1)$ (or more generally, of values $L(\chi_t,\alpha)$ with $\alpha\equiv 1\pmod{q-1}$). Already, we notice by Corollary \ref{fq} that
if $t$ belongs to $\FF_q^{\text{alg}}$, then $L(\chi_t,1)$ is transcendental by the well known transcendency of
$\widetilde{\pi}$. Moreover, the use of the functions 
$$f_{t}(X)=\prod_{i>0}(1-t/X^{q^i}),\quad |X|>1,$$ 
satisfying the functional equation $f_t(X)=(1-t/X^{q})f_t(X^q)$, the identity 
$L(\chi_t,1)=f_t(\theta)/f_\theta(\theta)$ and Mahler's method for transcendence as
in \cite{mahler}, allow to show that for all $t$ algebraic non-zero not of the form $\theta^{q^i}$ for $i\in\ZZ$, $L(\chi_t,1)$ and
$\widetilde{\pi}$ are algebraically independent. A precise necessary and
sufficient condition on $t_1,\ldots,t_s$ algebraic for the algebraic independence of
$\bsb{\Omega}(t_1),\ldots,\bsb{\Omega}(t_s)$ can be given, from which one deduces the 
corresponding condition for 
$L(\chi_{t_1},1),\ldots,L(\chi_{t_s},1)$. These facts will be described in more detail in another work.

\medskip

The proofs of Theorems \ref{corollairezeta11} and \ref{generalalpha} that we propose rely on certain properties
of {\em deformations of vectorial modular forms} (see Section \ref{dvmf}). In fact, Theorem \ref{corollairezeta11} is a corollary of 
an identity involving such functions that we describe now  (Theorem \ref{firsttheorem} below) and Theorem 
\ref{generalalpha} will be obtained from a simple modification of the techniques introduced to prove Theorem
\ref{firsttheorem}.

\medskip 

\noindent\emph{A fundamental identity for deformations of vectorial modular forms.} To present Theorem \ref{firsttheorem}, we need to introduce more tools. Let $\Omega$ be the rigid analytic
space $\CC_\infty\setminus K_\infty$. 
For $z\in\Omega$, we denote by $\Lambda_z$ the $A$-module $A+zA$, free of rank $2$.
The evaluation at $\zeta\in \CC_\infty$ of the {\em exponential function} $e_{\Lambda_z}$ associated to the lattice $\Lambda_z$ is given by the series 
\begin{equation}\label{ellipticexponential}
e_{\Lambda_z}(\eta)=\sum_{i=0}^\infty\alpha_i(z)\eta^{q^i},
\end{equation}
for  functions $\alpha_i:\Omega\rightarrow \CC_\infty$ with $\alpha_0=1$. We recall that for $i>0$, $\alpha_i$ is a 
{\em Drinfeld modular form} of weight $q^i-1$ and type $0$ in the sense of Gekeler, \cite{Ge}.

We also recall from \cite{archiv} the series:
\begin{eqnarray*}
\bsb{s}_1(z,t)&=&\sum_{i=0}^\infty\frac{\alpha_i(z)z^{q^i}}{\theta^{q^i}-t},\\
\bsb{s}_2(z,t)&=&\sum_{i=0}^\infty\frac{\alpha_i(z)}{\theta^{q^{i}}-t},
\end{eqnarray*} which converge for $|t|<q$ and define two
functions $\Omega\rightarrow \CC_\infty[[t]]$ with the series in the image converging for 
$|t|<q$. We point out that for a fixed choice of
$z\in\Omega$,
the matrix function ${}^t(\bsb{s}_1(z,t),\bsb{s}_2(z,t))$ is the {\em canonical rigid analytic trivialisation} of the {\em $t$-motive associated}
to the lattice $\Lambda_z$ discussed in \cite{Bourbaki}.
We set, for $i=1,2$:
$$\bsb{d}_i(z,t):=\widetilde{\pi}s_\carlitz(t)^{-1}\bfs_i(z,t),$$
remembering that in the notations of \cite{archiv}, we have $\bsb{d}_2=\bsb{d}$.
The advantage of using these functions instead of the $\bsb{s}_i$'s, is that 
$\bsb{d}_2$ has a {\em $u$-expansion} defined over $\FF_q[t,\theta]$ (see Proposition \ref{prarchiv}); moreover, 
it can 
be proved that for all $z\in\Omega$, $\bsb{d}_1(z,t)$ and $\bsb{d}_2(z,t)$ are entire functions
of the variable $t\in\CC_\infty$.

On the other hand, both the series 
$$\bsb{e}_1(z,t)=\sideset{}{'}\sum_{c,d\in A}\frac{\chi_t(c)}{cz+d},\quad \bsb{e}_2(z,t)=\sideset{}{'}\sum_{c,d\in A}\frac{\chi_t(d)}{cz+d}$$
converge for $(z,t)\in\Omega\times \CC_\infty$ with $|t|\leq 1$, where the dash $'$ denotes a sum avoiding the pair $(c,d)=(0,0)$. The series $\bsb{e}_1,\bsb{e}_2$ then define
functions $\Omega\rightarrow \CC_\infty[[t]]$ such that all the series in the images converge for $|t|\leq 1$.

Let $\mathbf{Hol}(\Omega)$ be the ring of holomorphic functions $\Omega\rightarrow\CC_\infty$, over which the Frobenius $\FF_q$-linear map $\tau$ is well defined: if $f\in\mathbf{Hol}(\Omega)$, then $\tau f=f^q$. We consider the
unique $\FF_q((t))$-linear extension of $\tau$: 
$$\mathbf{Hol}(\Omega)\otimes_{\FF_q}\FF_q((t))\rightarrow\mathbf{Hol}(\Omega)\otimes_{\FF_q}\FF_q((t)),$$
again denoted by $\tau$. 

We shall prove the fundamental theorem:
\begin{Theorem}\label{firsttheorem}
The following identities hold for $z\in\Omega,t\in\CC_\infty$ such that $|t|\leq 1$:
\begin{eqnarray}\label{eqP101}
L(\chi_t,1)^{-1}\bsb{e}_1(z,t)&=&-(t-\theta)s_\carlitz(t)(\tau\bsb{d}_2)(z,t) h(z),\\
L(\chi_t,1)^{-1}\bsb{e}_2(z,t)&=&(t-\theta)s_\carlitz(t)(\tau\bsb{d}_1)(z,t) h(z)\nonumber.\end{eqnarray}
\end{Theorem}
In the statement of the theorem, $h$ is the opposite of the unique normalised Drinfeld cusp form of weight $q+1$ and type $1$ for $\Gamma=\mathbf{GL}_2(A)$
as in \cite{Ge}. It can be observed that both right-hand sides in (\ref{eqP101}) are well defined for $t\in\CC_\infty$ with $|t|<q^q$ so that
these equalities provide analytic continuations of the functions $\bsb{e}_1,\bsb{e}_2$ in terms of
the variable $t$.

Let us write $\mathcal{E}=L(\chi_t,1)^{-1}(\bsb{e}_1,\bsb{e}_2)$ and $\mathcal{F}=\binom{\bsb{d}_1}{\bsb{d}_2}$, and
let us consider the representation $\rho_t:\mathbf{GL}_2(A)\rightarrow\mathbf{GL}_2(\FF_q[t])$ defined by
$$\rho_t(\gamma)=\sqm{\chi_t(a)}{\chi_t(b)}{\chi_t(c)}{\chi_t(d)},$$ for $\gamma=\sqm{a}{b}{c}{d}\in\mathbf{GL}_2(A)$.
Then, for any such a choice of $\gamma$ we have the functional equations:
\begin{eqnarray*}\mathcal{F}(\gamma(z))&=&(cz+d)^{-1}\rho_t(\gamma)\cdot\mathcal{F}(z),\\
{}^t\mathcal{E}(\gamma(z))&=&(cz+d)\;{}^t\rho_t^{-1}(\gamma)\cdot{}^t\mathcal{E}(z).\end{eqnarray*}
This puts right away the functions ${}^t\mathcal{E}$ and $\mathcal{F}$ in the framework of 
{\em deformations of vectorial modular forms}, topic that will be developed in Section \ref{dvmf}.

Let $B_1$ be the set of $t\in\CC_\infty$ such that $|t|\leq1$. 
We will make use of a remarkable sequence $\mathcal{G}=(\mathcal{G}_k)_{k\in\ZZ}$ of functions
$\Omega\times B_1\rightarrow\CC_\infty$ defined by the scalar 
product (with component-wise action of $\tau$):
$$\mathcal{G}_k=\mathcal{G}_{1,0,k}=(\tau^k\mathcal{E})\cdot\mathcal{F},$$ 
such that, for $k\geq 0$, $\mathcal{G}_k$ (\footnote{We will also adopt the notation $\mathcal{E}_{1,0}=\mathcal{E}$ and $\mathcal{G}_k=\mathcal{G}_{1,0,k}$.}) turns out to be an element of $M_{q^k-1,0}\otimes\FF_q[t,\theta]$ where $M_{w,m}$ denotes the $\CC_\infty$-vector space of Drinfeld modular forms 
of weight $w$ and type $m$. In fact, only the terms $\mathcal{G}_0,\mathcal{G}_1$ are needed to be examined to prove 
our Theorem. Once their explicit computation is accomplished (see Proposition \ref{interpretationvk}),
the proof of Theorem \ref{firsttheorem} is only a matter of solving a non-homogeneous system in
two equations and two indeterminates $\bsb{e}_1,\bsb{e}_2$. Furthermore, Theorem \ref{corollairezeta11}
will be deduced with the computation of a limit in the identity $\bsb{e}_1=-L(\chi_t,1)\tau (s_\carlitz\bsb{d}_2) h$
and a similar path will be followed for Theorem \ref{generalalpha}, by using this time the more general
sequences $\mathcal{G}_{\alpha,0,k}=(\tau^k\mathcal{E}_{\alpha,0})\cdot\mathcal{F}$ defined later.

\medskip

We end this introduction by pointing out that the sequence 
$\mathcal{G}$ has itself several interesting features. For example, the functions $\mathcal{G}_k$ already appeared 
in \cite{archiv} (they are denoted by $g_k^\star$ there) as the coefficients
of the ``cocycle terms $\bsb{L}_\gamma$" of the functional equations of the {\em deformations of Drinfeld quasi-modular forms $\tau^k\bsb{E}$} introduced there. 

It is also interesting to notice that the {\em deformation of Legendre's identity} (\ref{dethPsi}) that we quote here 
(proved in \cite{archiv}):
$$h^{-1}(\tau s_\carlitz)^{-1}=\bsb{d}_1(\tau\bsb{d}_2)-\bsb{d}_2(\tau\bsb{d}_1)$$
can be deduced from Theorem \ref{firsttheorem} by using the fact that $\mathcal{G}_0$ equals $-1$, property  obtained in 
Proposition \ref{interpretationvk}.

Moreover, it can be proved, with the help of the theory of {\em $\tau$-linear recurrent sequences} and {\em $\tau$-linearised recurrent sequences} (they will not be described here), that for $k\geq 0$, the function $\mathcal{G}_k(z,\theta)$, well defined, is equal to the {\em ortho-Eisenstein} series $g_k(z)$, and that $\mathcal{G}_k(z,\theta^{q^k})$, also well defined, is equal to the {\em para-Eisenstein} series $m_k(z)$, in the 
notations and the terminology of \cite{Gek2}. Hence, the sequence $\mathcal{G}$ provides an interesting tool also in the study of both these kinds of functions. This program will be however pursued in another paper.

\medskip

\noindent\emph{Acknowledgements.} The author is thankful to Vincent Bosser for a careful reading of a previous version of the manuscript, and Vincent Bosser, David Goss and Matt Papanikolas 
for fruitful discussions about the topics of the present paper.
 
\section{Vectorial modular forms and their deformations\label{dvmf}}

We denote by $J_\gamma$ the factor of automorphy $(\gamma,z)\mapsto cz+d$, if $\gamma=\sqm{a}{b}{c}{d}$. 
We will write, for $z\in\Omega$, $u=u(z)=1/e_\carlitz(\widetilde{\pi}z)$; this is the local 
parameter at infinity of $\Gamma\backslash\Omega$. For all $w,m$, we have an embedding of
$M_{w,m}$ in $\CC_\infty[[u]]$ associating to every Drinfeld modular form its $u$-expansion, see \cite{Ge};
we will often identify modular forms with their $u$-expansions.

In all the following, $t$ will be considered either as a new indeterminate or as a parameter varying in $\CC_\infty$, and we will freely switch from formal series to functions. 
We will say that a series $\sum_{i\geq i_0}c_iu^i$ (with the coefficients $c_i$ in some field extension of $\FF_q(t,\theta)$) is {\em normalised}, if $c_{i_0}=1$.
We will also say that the series is {\em of type} $m\in\ZZ/(q-1)\ZZ$ if $i\not\equiv m\pmod{q-1}$ implies $c_i=0$.
This definition is obviously compatible with the notion of {\em type} of a Drinfeld modular form already mentioned in the introduction, see \cite{Ge}.

Vectorial modular forms is a classical topic of investigation in the theory of modular forms of one complex 
variable. The following definition is a simple adaptation of the notion of vectorial modular form for $\mathbf{SL}_2(\ZZ)$ investigated in works by Knopp and Mason such as \cite{KM,Mas}, the most relevant for us.

Let $\rho:\Gamma\rightarrow\mathbf{GL}_s(\CC_\infty)$ be a representation of $\Gamma$.

\begin{Definition}
{\em A {\em vectorial modular form of weight $w$, type $m$ and dimension $s$ associated to $\rho$} is a
rigid holomorphic function $f:\Omega\rightarrow\mathbf{Mat}_{s\times 1}(\CC_\infty)$ with the following two properties. Firstly, for all
$\gamma\in\Gamma$,
$$f(\gamma(z))=\det(\gamma)^{-m}J_\gamma^w\rho(\gamma)\cdot f(z).$$ Secondly, the 
vectorial function $f={}^t(f_1,\ldots,f_s)$ is {\em tempered} at infinity in the following way. There exists
$\nu\in\ZZ$ such that for all $i\in\{1,\ldots,s\}$, the following limit holds: $$\lim_{|z|=|z|_i\rightarrow\infty}u(z)^\nu f_i(z)=0$$ ($|z|_i$ denotes, for $z\in \CC_\infty$, 
the infimum $\inf_{a\in K_\infty}\{|z-a|\}$).}
\end{Definition}
We denote by $\mathcal{M}^s_{w,m}(\rho)$ the $\CC_\infty$-vector space generated by these functions.
We further write, for $s=1$, $M^!_{w,m}=\mathcal{M}^1_{w,m}(\bsb{1})$ with $\rho=\bsb{1}$ the constant representation. This is just the $\CC_\infty$-vector space (of infinite dimension)
generated by the rigid holomorphic functions $f:\Omega\rightarrow\CC_\infty$ satisfying, for all $z\in\Omega$ and
$\gamma\in\Gamma$,
$f(\gamma(z))=\det(\gamma)^{-m}J_\gamma^wf(z)$ and meromorphic at infinity.
It can be proved that this vector space is generated by the functions $h^{-i}m_i$, where
$m_i$ is a Drinfeld modular form of weight $w+i(q+1)$ and type $i$. 


Next, we briefly discuss simple formulas expressing Eisenstein series for $\Gamma$ as scalar products
of two vectorial modular forms (\footnote{Similar identities hold in the classical theory of 
$\mathbf{SL}_2(\ZZ)$ but will not be discussed here.}); these formulas suggest the theory we develop here. 
Let $g_k$ be the normalisation of the Eisenstein series of weight $q^k-1$
as in \cite{Ge}. We have, for $k>0$ (see loc. cit.):
$$
g_k(z)=(-1)^{k+1}\widetilde{\pi}^{1-q^k}[k]\cdots[1]\left(z\sideset{}{'}\sum_{c,d\in A}\frac{c}{(cz+d)^{q^k}}+\sideset{}{'}\sum_{c,d\in A}\frac{d}{(cz+d)^{q^k}}\right).
$$
The identity above can then be rewritten in a more compact form as a scalar product:
\begin{equation}\label{makeup}g_k(z)=-\zeta(q^k-1)^{-1}\mathcal{E}_k'(z)\cdot\mathcal{F}'(z),\end{equation}
where $\mathcal{E}_k'(z)$ is the convergent
series $\sum_{c,d\in A}'(cz+d)^{-q^k}(c,d)$, defining a holomorphic map $\Omega\rightarrow \mathbf{Mat}_{1\times 2}(\CC_\infty)$  (it is well defined also
when $k=0$, but only as a conditionally convergent series), and
$\mathcal{F}'(z)$ is the map $\binom{z}{1}:\Omega\rightarrow \mathbf{Mat}_{2\times 1}(\CC_\infty)$ (\footnote{After this discussion we will no longer use the notations $\mathcal{E}_k'$ and $\mathcal{F}'$.}).
It is not difficult to see that ${}^t\mathcal{E}_k'$ and $\mathcal{F}'$ are {vectorial modular forms}.
More precisely, if $\rho$ is the identity map, we have ${}^t\mathcal{E}_k'\in\mathcal{M}^2_{q^k,0}({}^t\rho^{-1})$
and $\mathcal{F}'\in\mathcal{M}^2_{-1,0}(\rho)$ (\footnote{More generally, one speaks about matrix modular forms 
associated to left and right actions of two representations of $\Gamma$ (or $\mathbf{SL}_2(\ZZ)$). Then,
$\mathcal{E}_k'$ is a row matrix modular form associated to the right action of ${}^t\rho^{-1}$.}).

We are going to provide variants of the identities
(\ref{makeup}) depending on the variable $t\in\CC_\infty$ in such a way that they will become arithmetically rich. For this, we need
to introduce the notion of {\em deformation of vectorial modular form}.
In the rest of the paper, we will only use deformations of vectorial modular of a certain type.
For the sake of simplicity, we will confine the presentation to the class of functions that we strictly need.
This does not immediately allow to see how vectorial modular forms arise as special values of such deformations (for specific choices of 
$t\in\CC_\infty$). 
However, the reader should be aware that this, in an obvious extension of the theory presented
here, is of course possible.

The subring 
of formal series in $\CC_\infty[[t]]$ converging for all $t$ such that $|t|\leq 1$ will be denoted by $\TT$.
It is  endowed with the $\FF_q[t]$-linear automorphism $\tau$ acting on formal series as follows:
$$\tau\sum_ic_it^i=\sum_ic_i^qt^i.$$
We will work with certain functions $f:\Omega\times B_1\rightarrow \CC_\infty$ with 
the property that for all $z\in\Omega$,
$f(z,t)$ can be identified with a series of $\TT$ converging for all $t_0\in B_1$ to $f(z,t_0)$. For such functions we will then also write
$f(z)$ to stress the dependence on $z\in\Omega$ when we want to consider them as functions $\Omega\rightarrow\TT$.
Sometimes, we will not specify the variables $z,t$ and just write $f$ instead of $f(z,t)$ or $f(z)$ to lighten our formulas
just as we did in some parts of the introduction. Moreover,
$z$ will denote a variable in $\Omega$ all along the paper.

With $\mathbf{Hol}(\Omega)$ we denote the ring of rigid holomorphic functions on $\Omega$.
Let us denote by 
$\mathcal{R}$ the integral ring whose elements are the formal series $f=\sum_{i\geq0}f_it^i$,
such that
\begin{enumerate}
\item For all $i$, $f_i$ is a map $\Omega\rightarrow \CC_\infty$ belonging to $\mathbf{Hol}(\Omega)$.
\item For all $z\in\Omega$, $\sum_{i\geq0}f_i(z)t^i$ is an element of $\TT$.
\end{enumerate}
The ring $\mathcal{R}$ is
endowed with the injective endomorphism $\tau$ acting on formal series as follows:
$$\tau\sum_{i\geq0}f_i(z)t^i=\sum_{i\geq0}f_i(z)^qt^i.$$

\subsection{Deformations of vectorial modular forms.}



\medskip

Let us consider a representation
\begin{equation}\label{rho}
\rho:\Gamma\rightarrow \mathbf{GL}_s(\FF_q((t))).
\end{equation}
We assume that the determinant representation $\det(\rho)$ is the 
$\mu$-th power of the determinant character, for some $\mu\in\ZZ/(q-1)\ZZ$. 
\begin{Definition}
{\em A {\em deformation of vectorial modular form} of weight $w$, dimension $s$ and type $m$ associated with a representation $\rho$ as in 
(\ref{rho}) is a column matrix $\mathcal{F}={}^t(f_1,\ldots,f_s)\in\mathbf{Mat}_{s\times 1}(\mathcal{R})$ such that
the following two properties hold. Firstly, considering $\mathcal{F}$ as a map $\Omega\rightarrow\mathbf{Mat}_{s\times 1}(\TT)$ we have,
for all $\gamma\in\Gamma$, 
$$\mathcal{F}(\gamma(z))=J_\gamma^w\det(\gamma)^{-m}\rho(\gamma)\cdot\mathcal{F}(z).$$
Secondly, the entries of $\mathcal{F}$ are {\em tempered}: there exists $\nu\in\ZZ$ such that, for all $t\in B_1$ and $i\in\{1,\ldots,s\}$,
$\lim_{|z|=|z|_i\rightarrow\infty}u(z)^\nu f_i(z)=0$.}
\end{Definition}
The set of deformations of vectorial modular forms of weight $w$, dimension $s$ and type $m$ associated to a representation $\rho$
is a $\TT$-module that we will denote by $\mathcal{M}^{s}_{w,m}(\rho)$ (we use the same notations as
for vectorial modular forms). 
If $s=1$ and if $\rho=\bsb{1}$ is the constant map, then $\mathcal{M}^1_{w,m}(\bsb{1})$ is the space $M^!_{w,m}\otimes\TT$. It is easy to see that we can endow the space
$\mathcal{M}^2(\rho)=\oplus_{w,m}\mathcal{M}^2_{w,m}(\rho)$ with the structure of a graded
$M^!\otimes\TT$-module, where $M^!=\oplus_{w',m'}M^!_{w',m'}$.

\begin{Lemme}\label{twist} Let $k$ be a non-negative integer.
If $\mathcal{F}$ is in $\mathcal{M}^{s}_{w,m}(\rho)$, then $\tau^k\mathcal{F}\in\mathcal{M}^{s}_{wq^k,m}(\rho)$.
If we
choose nonnegative integers $k_1,\ldots,k_s$, then
$$\det(\tau^{k_1}\mathcal{F},\ldots,\tau^{k_s}\mathcal{F})\in M^!_{w(q^{k_1}+\cdots+q^{k_s}),sm-\mu}\otimes\TT.$$
In particular, $$W_\tau(\mathcal{F})=\det(\tau^0\mathcal{F},\ldots,\tau^{s-1}\mathcal{F})\in M^!_{w(1+q+q^2+\cdots+q^{s-1}),sm-\mu}\otimes\TT.$$
\end{Lemme}

\noindent\emph{Proof.} From the definition, and for all $k'\in\ZZ$,
$$(\tau^{k'}\mathcal{F})(\gamma(z))=J_\gamma^{wq^{k'}}\det(\gamma)^{-m}\rho(\gamma)(\tau^{k'}\mathcal{F})(z)$$
because $\tau(\rho(\gamma))=\rho(\gamma)$. Moreover, $\tau$ is an endomorphism of $\mathcal{R}$ and it
is obvious that $\mathcal{F}$ being tempered, also $\tau^k\mathcal{F}$ is tempered. 

Now define the matrix function:
$$\mathbf{M}_{k_1,\ldots,k_s}=(\tau^{k_1}\mathcal{F},\ldots,\tau^{k_s}\mathcal{F}).$$
After the first part of the lemma we have, for $\gamma\in\mathbf{GL}_2(A)$:
$$\mathbf{M}_{k_1,\ldots,k_s}(\gamma(z))=\det(\gamma)^{-m}\rho(\gamma)\cdot \mathbf{M}_{k_1,\ldots,k_s}(z)\cdot\mathbf{Diag}(J_\gamma^{wq^{k_1}},\cdots,
J_\gamma^{wq^{k_s}}),$$
and we conclude the proof taking determinants of both sides. \CVD

\begin{Lemme}\label{wronsk} Let 
us consider $\mathcal{F}$ in $\mathcal{M}^{s}_{w,m}(\rho)$ and let $\mathcal{E}$ be such that ${}^t\mathcal{E}$ is in $\mathcal{M}^{s}_{w',m'}({}^t\rho^{-1})$. Let us denote by $\mathcal{G}_k$ the scalar product $(\tau^k\mathcal{E})\cdot\mathcal{F}$, then, for nonnegative $k$, $$\mathcal{G}_k\in M^!_{wq^k+w',m+m'}\otimes\TT.$$ 
Furthermore, we have:
$$\tau^k\mathcal{G}_{-k}\in M^!_{w+w'q^k,m+m'}\otimes\TT.$$
\end{Lemme}
\noindent\emph{Proof.} By  Lemma \ref{twist},
$\tau^k({}^t\mathcal{E})$ is in $\mathcal{M}^s_{w'q^k,m'}({}^t\rho^{-1})$ and 
$\tau^k\mathcal{F}$ is in $\mathcal{M}^s_{wq^k,m}(\rho)$. 
Let 
$\gamma$ be in $\mathbf{GL}_2(A)$. We have, after transposition, and for all $k\in\ZZ$,
$$(\tau^{k}\mathcal{E})(\gamma(z))=J_\gamma^{w'q^{k}}\det(\gamma)^{-m'}\mathcal{E}(z)\cdot\rho^{-1}(\gamma),$$
and since $\tau^k\mathcal{G}_{-k}={}^t\mathcal{E}\cdot(\tau^k\mathcal{F})$,
$$(\tau^{k}\mathcal{F})(\gamma(z))=J_\gamma^{wq^{k}}\det(\gamma)^{-m}\rho(\gamma)\cdot(\tau^{k}\mathcal{F}(z)).$$
Hence, for $k\geq 0$,
$$\mathcal{G}_{k}(\gamma(z))=J_\gamma^{w'q^{k}+w}\det(\gamma)^{-m-m'}\mathcal{G}_{k}(z),$$
and
$$\tau^k\mathcal{G}_{-k}(\gamma(z))=J_\gamma^{w+w'q^{k}}\det(\gamma)^{-m-m'}\mathcal{G}_{-k}(z).$$
On the other hand, $\mathcal{G}_k$ and $\tau^k\mathcal{G}_{-k}$ are tempered for all $k\geq 0$, from which we can conclude.\CVD

From now on, we will use the representation $\rho=\rho_{t}$ and the transposed of its inverse. 

\subsection{The function $\mathcal{F}$}

The function of the title of this subsection is the vector valued function
$\binom{\bsb{d}_1}{\bsb{d}_2}$. In the next proposition, containing the properties of $\mathcal{F}$ of interest for us, we write $g$ for the unique normalised Drinfeld modular form of weight $q-1$ and type $0$ for $\Gamma$
(proportional to an Eisenstein series), and $\Delta$ for the cusp form $-h^{q-1}$.

\begin{Proposition}\label{prarchiv}
We have the following four properties for $\mathcal{F}$ and the $\bsb{d}_i$'s.

\begin{enumerate}
\item We have 
$\mathcal{F}\in\mathcal{M}^2_{-1,0}(\rho_{t})$.
\item The functions $\bsb{d}_1,\bsb{d}_2$ span the $\FF_q(t)$-vector space of dimension
$2$ of solutions of the following $\tau$-linear difference equation:
\begin{equation}\label{equexpansion}
X=(t-\theta^q)\Delta\tau^2X+g\tau X,
\end{equation}
in a suitable existentially closed inversive 
field containing $\mathcal{R}$.
\item Let us consider the matrix function:
$$\Psi (z,t):=\sqm{\bsb{d}_1(z,t)}{\bsb{d}_2(z,t)}{(\tau\bsb{d}_1)(z,t)}{(\tau\bsb{d}_2)(z,t)}.$$
For all $z\in\Omega$ and $t$ with $|t|<q$:
\begin{equation}\label{dethPsi}
\det(\Psi)=(t-\theta)^{-1}h(z)^{-1}s_\carlitz(t)^{-1}.
\end{equation}
\item We have the series expansion
\begin{equation}\label{uexpd2}
\bsb{d}_2 =\sum_{i\geq 0}c_i(t)u^{(q-1)i}\in1+u^{q-1}\FF_q[t,\theta][[u^{q-1}]],
\end{equation}
convergent for $t,u$ sufficiently close to $(0,0)$.
\end{enumerate}
\end{Proposition}
\noindent\emph{Proof.} All the properties but one follow immediately from the results of
\cite{archiv} where some of them are stated in slightly different, although equivalent formulations.
The only property we have to justify here is that $\mathcal{F}$ is tempered. After (\ref{uexpd2}),
we are led to check that there exists $\nu\in\ZZ$ such that $u(z)^\nu\bsb{d}_1\rightarrow0$
for $z\in\Omega$ such that $|z|=|z|_i\rightarrow\infty$. For this, we have the following lemma,
which concludes the proof of the Proposition.

\begin{Lemme}\label{basiclimits}
The following limits hold, for all $t\in\CC_\infty$ such that $|t|\leq 1$:
$$\lim_{|z|=|z|_i\rightarrow\infty}u(z)\bsb{d}_1(z,t)=0,\quad \lim_{|z|=|z|_i\rightarrow\infty}u(z)(\tau\bsb{d}_1)(z,t)=1.$$
\end{Lemme}
\noindent\emph{Proof.} We recall from \cite{archiv} the series expansion
$$\bsb{d}_1(z)=\frac{\widetilde{\pi}}{s_{\carlitz}(t)}\bsb{s}_2(z)=\frac{\widetilde{\pi}}{s_{\carlitz}(t)}\sum_{n\geq0}e_{\Lambda_z}\left(\frac{z}{\theta^{n+1}}\right)t^n,$$
converging for all $t$ such that $|t|<q$ and for all $z\in\Omega$.

By a simple modification of the proof of \cite[Lemma 5.9 p. 286]{gekeler:compositio}, we have
$$\lim_{|z|_i=|z|\to\infty}u(z)t^ne_{\Lambda_z}(z/\theta^{n+1})^q=0$$ uniformly in $n>0$, for all $t$ such that
$|t|\leq 1$ (in fact, even if $|t|\leq q$).

Moreover, it is easy to show that
$$\lim_{|z|_i=|z|\to\infty} u(z)e_{\Lambda_z}(z/\theta)^q=\widetilde{\pi}^{-q}\lim_{|z|_i=|z|\to\infty}e^{q}_{\carlitz}(\widetilde{\pi}z/\theta)/e_{\carlitz}(\widetilde{\pi}z)=1,$$ which gives the second limit, from which we deduce the first limit as well.\CVD

\subsection{Structure of $\mathcal{M}^2$}

Let us denote by $\mathcal{F}^*$ the function $\binom{-\bsb{d}_2}{\bsb{d}_1}$,
which is easily seen to be an element of $\mathcal{M}^2_{-1,-1}({}^t\rho_t^{-1})$.
In this subsection we give some information on the structure of the spaces 
$\mathcal{M}^2_{w,m}({}^t\rho_t^{-1})$.

\begin{Proposition}\label{embedding}
We have
$$\mathcal{M}^2_{w,m}({}^t\rho^{-1}_t)=(M^!_{w+1,m+1}\otimes\TT)\mathcal{F}^*\oplus
(M^!_{w+q,m+1}\otimes\TT)(\tau\mathcal{F}^*).$$ More precisely,
for all $\mathcal{E}$ with ${}^t\mathcal{E}\in\mathcal{M}^2_{w,m}({}^t\rho^{-1}_t)$ we have
$${}^t\mathcal{E}=(\tau s_\carlitz)h((\tau\mathcal{G}_{-1})\mathcal{F}^*+\mathcal{G}_{0}(\tau\mathcal{F}^*)),$$
where we have written $\mathcal{G}_k=(\tau^k({}^t\mathcal{E}))\cdot\mathcal{F}$ for all $k\in\ZZ$.
\end{Proposition}

The first part of the proposition is equivalent to the equality
$$\mathcal{M}^2_{w,m}(\rho_t)=(M^!_{w+1,m}\otimes\TT)\mathcal{F}\oplus
(M^!_{w+q,m}\otimes\TT)(\tau\mathcal{F}).$$ In the rest of this paper, we will only use the equality 
for $\mathcal{M}^2_{w,m}({}^t\rho^{-1}_t)$.

\medskip

\noindent\emph{Proof of Proposition \ref{embedding}.} Let us temporarily write $\mathcal{M}'$ for $(M^!_{w+1,m+1}\otimes\TT)\mathcal{F}^*+
(M^!_{w+q,m+1}\otimes\TT)(\tau\mathcal{F}^*)$. It is easy to show, thanks to the results of 
\cite{archiv}, that the sum is direct. Indeed, in loc. cit. it is proved that 
$\bsb{d}_2$, $\tau\bsb{d}_2$, $g$ and $h$ are algebraically independent over $\CC_\infty((t))$.
Moreover, $\mathcal{M}'$ clearly embeds in $\mathcal{M}^2_{w,m}({}^t\rho^{-1}_t)$. It remains 
to show the opposite inclusion.

By Proposition \ref{prarchiv}, the matrix $M=(\mathcal{F},\tau^{-1}\mathcal{F})$ is invertible.
From (\ref{dethPsi}) we deduce that
$$\tau M^{-1}=(t-\theta)s_\carlitz h\sqm{-\bsb{d}_2}{\bsb{d}_1}{\tau\bsb{d}_2}{-\tau\bsb{d}_1}.$$

Let $\mathcal{E}$ be such that ${}^t\mathcal{E}\in\mathcal{M}^2_{w,m}({}^t\rho^{-1}_t)$.
Thanks to the above expression for $\tau M^{-1}$, the identity:
$$\binom{\mathcal{E}}{\tau\mathcal{E}}\cdot\mathcal{F}=\binom{\mathcal{G}_0}{\mathcal{G}_1},$$ product of a $2\times 2$ matrix with a one-column matrix (which is the definition of $\mathcal{G}_0,\mathcal{G}_1$),
yields the formulas:
\begin{eqnarray}
\mathcal{E}& =&(\mathcal{G}_0,\tau^{-1}\mathcal{G}_1)\cdot M^{-1}\nonumber\\
&=&(t-\theta^{1/q})h^{1/q}(\tau^{-1}s_\carlitz)(\mathcal{G}_0,\tau^{-1}\mathcal{G}_1)\cdot
\left(\begin{array}{ll}-\tau^{-1}\bsb{d}_2&\tau^{-1}\bsb{d}_1\\ \bsb{d}_2&-\bsb{d}_1\end{array}\right)\nonumber\\
&=&(t-\theta^{1/q})h^{1/q}(\tau^{-1}s_\carlitz)(\tau^{-1}\mathcal{G}_1\bsb{d}_2  -\mathcal{G}_0\tau^{-1}\bsb{d}_2,-\tau^{-1}\mathcal{G}_1\bsb{d}_1  +\mathcal{G}_0\tau^{-1}\bsb{d}_1).\label{partial1}
\end{eqnarray}
Now, we observe that we have, from the second part of Proposition \ref{prarchiv} and for all $k\in\ZZ$, 
\begin{equation}\nonumber
\mathcal{G}_k=g\tau\mathcal{G}_{k-1}+\Delta(t-\theta^q)\tau^2\mathcal{G}_{k-2}.
\end{equation}
Applying this formula for $k=1$ we obtain
\begin{equation}\label{taudifference}
\tau\mathcal{G}_{-1}=\frac{\tau^{-1}\mathcal{G}_{1}-g^{1/q}\mathcal{G}_{0}}{(t-\theta)\Delta^{1/q}},
\end{equation}
 and by using again Part two of Proposition \ref{prarchiv}, we eliminate $\tau^{-1}\mathcal{G}_{1}$
 and $\tau^{-1}\bsb{d}_i$:
$$(\tau^{-1}\mathcal{G}_{1})\bsb{d}_i  -\mathcal{G}_{0}\tau^{-1}\bsb{d}_i=\Delta^{1/q}(t-\theta)((\tau\mathcal{G}_{-1})\bsb{d}_i-\mathcal{G}_{0}(\tau\bsb{d}_i)),\quad i=1,2.$$
Replacing in (\ref{partial1}), and using $\Delta^{1/q}h^{1/q}=-h$ and $(t-\theta^{1/q})\tau^{-1}s_{\carlitz}=s_\carlitz$, 
we get the formula:
\begin{equation}\label{formula1}
\mathcal{E}=(t-\theta)s_\carlitz h(-(\tau\mathcal{G}_{-1})\bsb{d}_2+\mathcal{G}_{0}(\tau\bsb{d}_2),(\tau\mathcal{G}_{-1})\bsb{d}_1-\mathcal{G}_{0}(\tau\bsb{d}_1)).
\end{equation}
By Lemma \ref{wronsk}, we have $\mathcal{G}_{0}\in M^!_{w-1,m}\otimes\TT,(\tau\mathcal{G}_{-1})\in M^!_{w-q,m}\otimes\TT$ and the proposition follows from the fact that $h\in M_{q+1,1}$.\CVD

\noindent\emph{Remark.} Let us choose any $\mathcal{E}=(f_1,f_2)$ such that 
${}^t\mathcal{E}\in\mathcal{M}^{2}_{w,m}({}^t\rho_t^{-1})$.
Proposition \ref{embedding} implies that there exists $\mu\in\ZZ$ such that
$$h^\mu f_1\in\MM^\dag_{\mu(q+1)+w,1,\mu+m},$$
where $\MM^\dag_{\alpha,\beta,m}$ is a sub-module of {\em almost $A$-quasi-modular forms}
as introduced in \cite{archiv}.


\subsection{Deformations of vectorial Eisenstein and Poincar\'e series\label{poincare}.}
The aim of this subsection is to construct non-trivial elements 
of $\mathcal{M}^2_{w,m}({}^t\rho_t^{-1})$.

Following Gekeler \cite[Section 3]{Ge}, we recall that for all $\alpha>0$ there exists a polynomial $G_\alpha(u)\in \CC_\infty[u]$, called the 
{\em $\alpha$-th Goss polynomial},
such that, for all $z\in\Omega$, $G_{\alpha}(u(z))$ equals the sum of the convergent series
$$\widetilde{\pi}^{-\alpha}\sum_{a\in A}\frac{1}{(z+a)^\alpha}.$$

Several properties of these polynomials are collected in \cite[Proposition (3.4)]{Ge}. Here, we will need that
for all $\alpha$, $G_\alpha$ is of type $\alpha$ as a formal series of $\CC_\infty[[u]]$. Namely:
$$G_\alpha(\lambda u)=\lambda^\alpha G_\alpha(u),\quad \text{ for all }\lambda\in\FF_q^\times.$$

We also recall, for $a\in A$, the function
$$u_a(z):=u(az)=e_\carlitz(\widetilde{\pi}az)^{-1}=u^{|a|}f_a(u)^{-1}=u^{|a|}+\cdots\in A[[u]],$$
where $f_a\in A[[u]]$ is the {\em $a$-th inverse cyclotomic polynomial} defined in \cite[(4.6)]{Ge}. Obviously,
we have
$$u_{\lambda a}=\lambda^{-1}u_a\quad \text{ for all }\lambda\in\FF_q^\times.$$

We will use the following lemma.

\begin{Lemme}\label{primolemma}
Let $\alpha$ be a positive integer such that $\alpha\equiv 1\pmod{q-1}$.
We have, for all $t\in \CC_\infty$ such that $|t|\leq1$ and $z\in\Omega$, convergence of the series below, and equality:
$$\sideset{}{'}\sum_{c,d\in A}\frac{\ol{c}}{(cz+d)^{\alpha}}=-\widetilde{\pi}^\alpha\sum_{c\in A^+}\ol{c}G_{\alpha}(u_c(z)),$$ from which it follows that the series in the left-hand side is not identically zero.
\end{Lemme}
\noindent\emph{Proof.} Convergence is ensured by Lemma \ref{interpret} (or Proposition \ref{mainproppoincare})
and the elementary properties of Goss' polynomials. By the way, the series on the right-hand side 
converges for all $t\in\CC_\infty$.
We then compute:
\begin{eqnarray*}
\sideset{}{'}\sum_{c,d}\frac{\ol{c}}{(cz+d)^\alpha}&=&\sum_{c\neq0}\ol{c}\sum_{d\in A}\frac{1}{(cz+d)^\alpha}\\
&=&\widetilde{\pi}^\alpha\sum_{c\neq 0}\ol{c}\sum_{d\in A}\frac{1}{(c\widetilde{\pi}z+d\widetilde{\pi})^\alpha}\\
&=&\widetilde{\pi}^\alpha\sum_{c\neq 0}\ol{c}G_\alpha(u_c)\\
&=&\widetilde{\pi}^\alpha\sum_{c\in A^+}\ol{c}G_\alpha(u_c)\sum_{\lambda\in\FF_q^\times}\lambda^{1-\alpha}\\
&=&-\widetilde{\pi}^\alpha\sum_{c\in A^+}\ol{c}G_\alpha(u_c).
\end{eqnarray*} The non-vanishing of the series comes from the non-vanishing contribution of the 
term $G_{\alpha}(u)$ in the last series. Indeed, the order of vanishing at $u=0$ of the right-hand side is equal to 
$\min_{c\in A^+}\{\text{ord}_{u=0}G_\alpha(u_c)\}=\text{ord}_{u=0}G_\alpha(u)$, which is 
$<\infty$.
\CVD

Following \cite{Ge}, we consider the subgroup $$H=\left\{\sqm{*}{*}{0}{1}\right\}$$ of $\Gamma=\mathbf{GL}_{2}(A)$ and its left action on $\Gamma$.
For $\delta=\sqm{a}{b}{c}{d}\in\Gamma$, the map $\delta\mapsto(c,d)$ induces a bijection
between the orbit set $H\backslash\Gamma$ and the set of $(c,d)\in A^2$ with $c,d$ relatively prime. 

We consider the factor of automorphy $$\mu_{\alpha,m}(\delta,z)=\det(\delta)^{-m}J_\delta^\alpha,$$ where $m$ and $\alpha$ are 
positive integers (at the same time, $m$ will also determine a type, that is, a class modulo $q-1$).

Let $V_1 (\delta)$ be the row matrix $(\ol{c},\ol{d}).$
It is easy to show that the row matrix $$\mu_{\alpha,m}(\delta,z)^{-1}u^m(\delta(z))V_1 (\delta)$$ only depends on the class of $\delta\in H\backslash\Gamma$,
so that we can consider the following expression:
$$\mathcal{E}_{\alpha,m}(z)=\sum_{\delta\in H\backslash\Gamma}\mu_{\alpha,m}(\delta,z)^{-1}u^m(\delta(z))V_1 (\delta),$$ 
which is a row matrix  whose two entries are formal series. 

Let $\mathcal{V}$ be the set of functions $\Omega\rightarrow\mathbf{Mat}_{1\times 2}(\CC_\infty[[t]]).$
We introduce, for $\alpha,m$ integers, $f\in \mathcal{V}$ and $\gamma\in\Gamma$, the Petersson slash operator:
$$f|_{\alpha,m}\gamma=\det(\gamma)^{m}(cz+d)^{-\alpha}f(\gamma(z))\cdot\rho_{t}(\gamma).$$
This will be used in the next proposition, where we denote by $\log^+_q(x)$ the maximum between $0$ and $\log_q(x)$, the logarithm in base $q$  of $x>0$. We point out that we will not apply this proposition in full generality.
Indeed, in this paper, we essentially consider the case $m=0$ in the proposition. The proposition is presented
in this way for the sake of completeness.

\begin{Proposition}\label{mainproppoincare} Let $\alpha,m$ be non-negative integers with $\alpha\geq 2m+1$, and write $r(\alpha,m)=\alpha-2m-1$. We have the following properties.

\begin{enumerate}

\item For $\gamma\in\Gamma$, the map $f\mapsto f|_{\alpha,m}\gamma$ induces a permutation of the subset of $\mathcal{V}$: 
$$\mathcal{S}=\{\mu_{\alpha,m}(\delta,z)^{-1}u^m(\delta(z))V_1 (\delta);\delta\in H\backslash\Gamma\}.$$

\item If $t\in \CC_\infty$ and $\alpha,m$ are chosen so that $r(\alpha,m)>\log^+_q|t|$, then the components of $\mathcal{E}_{\alpha,m}(z,t)$ are series of functions
of $z\in\Omega$ which converge absolutely and uniformly on every compact subset of $\Omega$ to holomorphic functions.

\item If $|t|\leq 1$, then the components of $\mathcal{E}_{\alpha,m}(z,t)$ converge absolutely and uniformly on every compact subset of $\Omega$ 
also if $\alpha-2m>0$.

\item For any choice of $\alpha,m,t$ submitted to the convergence conditions above, the function ${}^t\mathcal{E}_{\alpha,m}(z,t)$ 
belongs to the space $\mathcal{M}^{2}_{\alpha,m}({}^t\rho_{t}^{-1})$.

\item If $\alpha-1\not\equiv 2m\pmod{(q-1)}$, the matrix function $\mathcal{E}_{\alpha,m}(z,t)$ is identically zero.

\item If $\alpha-1\equiv 2m\pmod{(q-1)}$, $\alpha\geq (q+1)m+1$ so that $\mathcal{E}_{\alpha,m}$ converges, then $\mathcal{E}_{\alpha,m}$ is not identically zero in its domain of convergence.
\end{enumerate}
\end{Proposition}

\noindent\emph{Proof.} 
\noindent\emph{1.} We choose $\delta\in H\backslash \Gamma$ corresponding
to a pair $(c,d)\in A^2$ with $c,d$ relatively prime, and set $f_\delta=\mu_{\alpha,m}(\delta,z)^{-1}u^m(\delta(z))V_1 (\delta)\in\mathcal{S}$.
We have
\begin{eqnarray*}
f_\delta(\gamma(z))&=&\mu_{\alpha,m}(\delta,\gamma(z))^{-1}u^m(\delta(\gamma(z)))V_1 (\delta)\\
&=&\mu_{\alpha,m}(\gamma,z)\mu_{\alpha,m}(\delta\gamma,z)^{-1}u^m(\delta\gamma(z))V_1 (\delta),\\
&=&\mu_{\alpha,m}(\gamma,z)\mu_{\alpha,m}(\delta\gamma,z)^{-1}u^m(\delta\gamma(z))V_1 (\delta\gamma)\cdot\rho_{t}(\gamma)^{-1},\\
&=&\mu_{\alpha,m}(\gamma,z)\mu_{\alpha,m}(\delta',z)^{-1}u^m(\delta'(z))V_1 (\delta')\cdot\rho_{t}(\gamma)^{-1},\\
&=&\mu_{\alpha,m}(\gamma,z)f_{\delta'}\cdot\rho_{t}(\gamma)^{-1},
\end{eqnarray*}
with $\delta'=\delta\gamma$ and $f_{\delta'}=\mu_{\alpha,m}(\delta',z)^{-1}u^m(\delta'(z))V_1 (\delta')$, from which part 1 of the
proposition follows.

\medskip

\noindent\emph{2.}
Convergence and holomorphy are ensured by simple modifications of \cite[(5.5)]{Ge}, or by the arguments in \cite[Chapter 10]{GePu}.
More precisely, let us choose $0\leq s\leq 1$ and look at the component at the place $s+1$ $$\mathcal{E}_s(z,t)=\sum_{\delta\in H\backslash\Gamma}\mu_{\alpha,m}(\delta,z)^{-1}u(\delta(z))^m\ol{c^sd^{1-s}}$$ of the vector 
series $\mathcal{E}_{\alpha,m}$. Writing $\alpha=n(q-1)+2m+l'$ with $n$ non-negative integer and $l'\geq 1$ we see, following Gerritzen and van der Put, \cite[pp. 304-305]{GePu} and taking into account the inequality
$|u(\delta(z))|\leq|cz+d|^2/|z|_i$, that
the term 
of the series $\mathcal{E}_s$:
$$\mu_{\alpha,m}(\delta,z)^{-1}u^m(\delta(z))\ol{c^sd^{1-s}}=(cz+d)^{-n(q-1)-l'-2m}u(\delta(z))^m\chi_t(c^sd^{1-s})$$
 (where $\delta$ corresponds to $(c,d)$)
has absolute value bounded from above by 
$$|z|_i^{-m}\left|\frac{\ol{c^sd^{1-s}}}{(cz+d)^{n(q-1)+l'}}\right|.$$
Applying the first part of the proposition, 
to check convergence, we can freely substitute $z$ with $z+a$ with $a\in A$ and we may assume, without loss of generality, that $|z|=|z|_i$. 
We verify that, either $\lambda=\deg_\theta z\in\QQ\setminus\ZZ$,
or $\lambda\in\ZZ$ case in which for all $\rho\in K_\infty$ with $|\rho|=|z|$, we have $|z-\rho|=|z|$.
In both cases, for all $c,d$, $|cz+d|=\max\{|cz|,|d|\}$.
Then, the series defining $\mathcal{E}_s$ can be decomposed as follows:
$$\mathcal{E}_s=\left(\sideset{}{'}\sum_{|cz|<|d|}+\sideset{}{'}\sum_{|cz|\geq|d|}\right)\mu_{\alpha,m}(\delta,z)^{-1}u^m(\delta(z))\ol{c^sd^{1-s}}.$$
We now look for upper bounds for the absolute values of the terms of the series above separating the two cases in a way similar to that of Gerritzen and van der Put
in loc. cit. 

Assume first that $|cz|<|d|$, that is, $\deg_\theta c+\lambda< \deg_\theta d$. Then
$$\left|\frac{\ol{c^sd^{1-s}}}{(cz+d)^{n(q-1)+l'}}\right|\leq \kappa \max\{1,|t|\}^{\deg_\theta d}|d|^{-n(q-1)-l'}\leq \kappa q^{\deg_\theta d(\log^+_q|t|-n(q-1)-l')},$$
where $\kappa$ is a constant depending on $\lambda$, and the corresponding sub-series 
converges with the imposed conditions on the parameters, because $\log^+_q|t|-n(q-1)-l'<0$.

If on the other side $|cz|\geq|d|$, that is, $\deg_\theta c+\lambda\geq\deg_\theta d,$ then
$$\left|\frac{\ol{c^sd^{1-s}}}{(cz+d)^{n(q-1)+l'}}\right|\leq\kappa' \max\{1,|t|\}^{\deg_\theta d}|c|^{-n(q-1)-l'}\leq\kappa' q^{\deg_\theta c(\log^+_q|t|-n(q-1)-l')},$$
with a constant $\kappa'$ depending on $\lambda$, again because $\log^+_q|t|-n(q-1)-l'<0$. This completes the proof of the second part of the Proposition.

\medskip

\noindent\emph{3.} This property can be deduced from the proof of the second part because if $|t|\leq 1$, then $|\chi_t(c^sd^{1-s})|\leq 1$. 

\medskip

\noindent\emph{4.} The property is obvious by the first part of the proposition, because 
$\mathcal{E}_{\alpha,m}=\sum_{f\in\mathcal{S}}f$, and because the functions are obviously
tempered thanks to the estimates we used in the proof of Part two.

\medskip

\noindent\emph{5.}
We consider $\gamma=\mathbf{Diag}(1,\lambda)$ with $\lambda\in\FF_q^\times$; the 
corresponding homography, multiplication by $\lambda^{-1}$, is equal to that defined by $\mathbf{Diag}(\lambda^{-1},1)$. Hence, we have:
\begin{eqnarray*}
\mathcal{E}_{\alpha,m}(\gamma(z))&=&\lambda^{\alpha-m}\mathcal{E}_{\alpha,m}(z)\cdot\mathbf{Diag}(1,\lambda^{-1})\\
&=&\lambda^{m}\mathcal{E}_{\alpha,m}(z)\cdot\mathbf{Diag}(\lambda,1),
\end{eqnarray*}
from which it follows that $\mathcal{E}_{\alpha,m}$ is identically zero if $\alpha-1\not\equiv 2m\pmod{q-1}$.

\medskip

\noindent\emph{6.} If $m=0$ and $\alpha=1$, we simply appeal to Lemma \ref{primolemma}.
Assuming now that either $m>0$ or $\alpha>1$, we have that $\mathcal{E}_{\alpha,m}$ converges 
at $t=\theta$ and:
\begin{eqnarray*}
z\mathcal{E}_0(z,\theta)+\mathcal{E}_1(z,\theta)&=&\sum_{\delta\in H\backslash\Gamma}\det(\delta)^{m}(cz+d)^{1-\alpha}u(\delta(z))^m\\
&=&P_{\alpha-1,m},
\end{eqnarray*}
where $P_{\alpha-1,m}\in M_{\alpha-1,m}$ is the Poincar\'e series of weight $\alpha-1$ and type $m$ so that
\cite[Proposition 10.5.2]{GePu} suffices for our purposes.\CVD

Let $\alpha,m$ be non-negative integers such that $\alpha-2m>1$ and $\alpha-1\equiv 2m\pmod{(q-1)}$. We have constructed functions:
\begin{eqnarray*}
\mathcal{E}_{\alpha,m}:\Omega&\rightarrow&\mathbf{Mat}_{1\times 2}(\mathcal{R}),\\
\mathcal{F}:\Omega&\rightarrow&\mathbf{Mat}_{2\times 1}(\mathcal{R}),
\end{eqnarray*}
and 
${}^t\mathcal{E}_{\alpha,m}\in\mathcal{M}^{2}_{\alpha,m}({}^t\rho_{t}^{-1})$, $\mathcal{F}\in\mathcal{M}^{2}_{-1,0}(\rho_{t})$. Therefore,
after Lemma \ref{wronsk}, the functions $$\mathcal{G}_{\alpha,m,k}=(\tau^k\mathcal{E}_{\alpha,m})\cdot\mathcal{F}=\mathcal{E}_{q^k\alpha,m}\cdot\mathcal{F}:\Omega\rightarrow\TT$$ satisfy $\mathcal{G}_{\alpha,m,k}\in M^!_{q^k\alpha-1,m}\otimes\TT$.

\medskip

\noindent\emph{A special case.}
After Proposition \ref{mainproppoincare},
if $\alpha>0$ and $\alpha\equiv 1\pmod{q-1}$, then $\mathcal{E}_{\alpha,0}\neq0$. We call these series {\em deformations of vectorial Eisenstein series}.

\begin{Lemme}\label{interpret} With $\alpha>0$ such that $\alpha\equiv 1\pmod{q-1}$, the following identity holds
for all $t\in\CC_\infty$ such that $|t|\leq 1$:
$$\mathcal{E}_{\alpha,0}(z,t)=L(\chi_t,\alpha)^{-1}\sideset{}{'}\sum_{c,d}(cz+d)^{-\alpha} V_1(c,d),$$
and $\mathcal{E}_{\alpha,0}$ is not identically zero.
\end{Lemme}
\noindent\emph{Proof.}
We recall the notation $$V_1(c,d)=(\ol{c},\ol{d})\in\mathbf{Mat}_{1\times 2}(\FF_q[t]).$$ 
We have
\begin{eqnarray*}
\sideset{}{'}\sum_{c,d}(cz+d)^{-\alpha} V_1(c,d)&=& \sum_{(c',d')=1}\sum_{a\in A^+}a^{-\alpha}(c'z+d')^{-\alpha} V_1(ac',ad')\\
&=&L(\chi_t,\alpha)\mathcal{E}_{\alpha,0}(z,t),
\end{eqnarray*}
where the first sum is over pairs of $A^2$ distinct from $(0,0)$, while the second sum is over the pairs $(c',d')$ of relatively prime
elements of $A^2$.
Convergence features are easy to deduce from Proposition \ref{mainproppoincare}. 
Indeed, we have convergence if $\log_q^+|t|<r(\alpha,m)=\alpha-1$, that is, $\max\{1,|t|\}\leq q^{\alpha-1}$ if $\alpha>1$
and we have convergence, for $\alpha=1$, for $|t|\leq1$. In all cases, convergence holds for $|t|\leq1$. 
Non-vanishing of the function also follows from Proposition \ref{mainproppoincare}.
\CVD


\section{Proof of the Theorems\label{proofs}}

We will need two auxiliary lemmas.

\begin{Lemme}\label{lemmelimit1} Let $\alpha>0$ be an integer such that $\alpha\equiv1\pmod{q-1}$. For all $t\in \CC_\infty$ such that $|t|\leq1$, we have
$$\lim_{|z|_i=|z|\to\infty}\bsb{d}_1(z)\sideset{}{'}\sum_{c,d}\frac{\ol{c}}{(cz+d)^\alpha}=0.$$
\end{Lemme}
\noindent\emph{Proof.} By Lemma \ref{basiclimits}, we have that 
$\lim_{|z|=|z|_i\rightarrow\infty}f(z)\bsb{d}_1(z,t)=0$ for all $t\in B_1$ and for all $f$ of the form
$f(z)=\sum_{n=1}^\infty c_nu(z)^n$, with $c_i\in\CC_\infty$, locally convergent at $u=0$.
But after Lemma \ref{primolemma}, $\sum_{c,d}'\frac{\ol{c}}{(cz+d)^\alpha}$ is equal, for
$|z|_i$ big enough, to the sum of the series 
$f(z)=-\widetilde{\pi}^\alpha\sum_{c\in A^+}\ol{c}G_{\alpha}(u_c(z))$ which is of the 
form $-\widetilde{\pi}^\alpha u^{\alpha}+o(u^{\alpha})$,
 and the lemma follows.\CVD 

\begin{Lemme}\label{lemmelimit2} Let $\alpha>0$ be an integer such that $\alpha\equiv1\pmod{q-1}$. For all $t\in \CC_\infty$ such that $|t|\leq1$, we have
$$\lim_{|z|_i=|z|\to\infty}\sideset{}{'}\sum_{c,d}\frac{\ol{d}}{(cz+d)^\alpha}=-L(\chi_t,\alpha).$$
\end{Lemme}
\noindent\emph{Proof.} It suffices to 
show that
$$\lim_{|z|_i=|z|\to\infty}\sum_{c\neq0}\sum_{d\in A}\frac{\ol{d}}{(cz+d)^{\alpha}}=0.$$
Assuming that $z'\in\Omega$ is such that $|z'|=|z'|_i$, we see that for all $d\in A$, $|z'+d|\geq |z'|_i$.
Now, consider $c\in A\setminus\{0\}$ and $z'=cz$ with $|z|=|z|_i$. Then,
$|cz+d|\geq|cz|_i=|cz|$, so that, for $|t|\leq 1$,
$$
\left|\frac{\chi_t(d)}{(cz+d)^{\alpha}}\right|\leq|cz|^{-\alpha}.$$
This implies that 
$$\left|\sum_{c\neq0}\sum_{d\in A}\frac{\ol{d}}{(cz+d)^{\alpha}}\right|\leq|z|^{-\alpha},$$
from which the Lemma follows.\CVD

The next step is to prove the following proposition.

\begin{Proposition}\label{interpretationvk}
For all $\alpha>0$ with $\alpha\equiv 1\pmod{q-1}$, then $\mathcal{G}_{\alpha,0,0}\in M_{\alpha-1,0}\otimes\TT$, and we have the limit $\lim_{|z|=|z|_i\rightarrow\infty}\mathcal{G}_{\alpha,0,0}=-1$.

Moreover, if $\alpha\leq q(q-1)$,
then:
$$\mathcal{G}_{\alpha,0,0}=-E_{\alpha-1},$$
where $E_{\alpha-1}$ is the normalised Eisenstein series of weight $\alpha-1$ for $\Gamma$.
\end{Proposition}

\noindent\emph{Proof.} Let us write:
$$F_\alpha(z,t):=\bsb{d}_1(z)\sideset{}{'}\sum_{c,d}\frac{\ol{c}}{(cz+d)^\alpha}+\bsb{d}_2(z)\sideset{}{'}\sum_{c,d}\frac{\ol{d}}{(cz+d)^\alpha},$$
series that converges for all $(z,t)\in\Omega\times\CC_\infty$ with $|t|\leq 1$. 
By Lemma \ref{interpret}, we have
$$F_\alpha(z,t)=L(\chi_t,\alpha)\mathcal{E}_{\alpha,0}(z,t)\cdot\mathcal{F}(z,t)=L(\chi_t,\alpha)\mathcal{G}_{\alpha,0,0},$$ so that $F_\alpha\in M^!_{\alpha-1,0}\otimes\TT$.
After (\ref{uexpd2}), we verify that for all $t$ with $|t|\leq1$,
$\lim_{|z|_i=|z|\to\infty}\bsb{d}_2(z)=1$. From Lemmas \ref{lemmelimit1} and \ref{lemmelimit2},
$$\lim_{|z|_i=|z|\to\infty}F_\alpha(z,t)=-L(\chi_t,\alpha).$$

Therefore, for all $t$ such that $|t|\leq 1$, $F_\alpha(z,t)$ converges to an holomorphic function
on $\Omega$ and is endowed with a $u$-expansion holomorphic at infinity.
In particular, $F_\alpha(z,t)$ is a family of modular forms of $M_{\alpha-1,0}\otimes\TT$.

Since for the selected values of $\alpha$, $M_{\alpha-1,0}=\langle E_{\alpha-1}\rangle$, 
we obtain that $F_\alpha=-L(\chi_t,\alpha)E_{\alpha-1}$.\CVD

\medskip

\noindent\emph{Proof of Theorem \ref{firsttheorem}.} 
Let us consider, for given $\alpha>0$, the form $\mathcal{E}=\mathcal{E}_{\alpha,0}$ and the scalar product form $\mathcal{G}_{\alpha,0,k}=(\tau^k\mathcal{E})\cdot\mathcal{F}$.
The general computation of $\mathcal{G}_0=\mathcal{G}_{\alpha,0,0}$
and $\tau\mathcal{G}_{-1}=\tau\mathcal{G}_{\alpha,0,-1}$ as in Proposition \ref{embedding}
is difficult, but for $\alpha=1$ we can apply 
Proposition \ref{interpretationvk}. We have $\mathcal{G}_{1,0,0}=-1$ and $\mathcal{G}_{1,0,1}=
\mathcal{G}_{q,0,0}=-g=-E_{q-1}$. Therefore, $\mathcal{G}_{\alpha,0,-1}=0$ by (\ref{taudifference}) and Theorem \ref{firsttheorem} follows.\CVD

\noindent\emph{Proof of Theorem \ref{corollairezeta11}.} 
Lemma \ref{primolemma} and Proposition \ref{embedding} imply that
\begin{equation}\label{Lchit}
L(\chi_t,\alpha)=\frac{\widetilde{\pi}^\alpha\sum_{c\in A^+}\ol{c}G_{\alpha}(u_c)}{(\tau s_\carlitz) h((\tau\mathcal{G}_{\alpha,0,-1})\bsb{d}_2-\mathcal{G}_{\alpha,0,0}(\tau \bsb{d}_2))},\end{equation} so that 
in particular, the numerator and the denominator of the fraction are one proportional to the other.
For $\alpha=1$ we can replace, by the above discussion, $\mathcal{G}_{\alpha,0,-1}=0$ and 
$\mathcal{G}_{\alpha,0,0}=-1$ and
we get, thanks to the fact that $h=-u+o(u)$ and $\sum_{c\in A^+}\chi_t(c)u_c=u+o(u)$,
$$L(\chi_t,\alpha)=\frac{\widetilde{\pi}\sum_{c\in A^+}\ol{c}u_c}{(\tau (s_\carlitz\bsb{d}_2)) h}=-\frac{\widetilde{\pi}}{\tau s_\carlitz}+o(1),$$
from which we deduce Theorem \ref{corollairezeta11} and even some additional information, namely, the formula:
$$(\tau \bsb{d}_2) h=-\sum_{c\in A^+}\ol{c}u_c.$$\CVD

\noindent\emph{Proof of Theorem \ref{generalalpha}.} For general $\alpha$, we know by (\ref{Lchit}) that there exists $\lambda_\alpha\in\LL$ ($\LL$ being the fraction field of $\TT$) such that
\begin{equation}\label{secondnewformula}
\sum_{c\in A^+}\ol{c}G_{\alpha}(u_c)=\lambda_\alpha h((\tau\mathcal{G}_{\alpha,0,-1})\bsb{d}_2-\mathcal{G}_{\alpha,0,0}(\tau \bsb{d}_2)).\end{equation} Since 
$L(\chi_t,\alpha)=\lambda_\alpha\widetilde{\pi}^\alpha/(\tau s_\carlitz)$, it remains to show that $\lambda_\alpha$
belongs to $\FF_q(t,\theta)$. 

Let us write $f$ for the series $\sum_{c\in A^+}\ol{c}G_{\alpha}(u_c)$,
$\phi$ for $\lambda_\alpha h\tau\mathcal{G}_{\alpha,0,-1}$ and $\psi$ for $-\lambda_\alpha h\mathcal{G}_{\alpha,0,0}$, so that 
(\ref{secondnewformula}) becomes:
$$f=\phi\bsb{d}_2+\psi\tau \bsb{d}_2.$$ Proposition \ref{embedding} then tells us that $\phi\in M^!_{\alpha+1,1}\otimes\LL$ and $\psi\in M^!_{\alpha+q,1}\otimes\LL$.

Let $L$ be an algebraically closed 
field containing $\LL$, hence containing also $\FF_q(t,\theta)$.
As for any choice of $w,m$, $M^!_{w,m}$ embeds in $\CC_\infty((u))$ and since there is a basis of this space with $u$-expansions defined over $K$, we have that $\mathbf{Aut}(L/\FF_q(t,\theta))$
acts on $M^!_{w,m}\otimes L$ through the coefficients of the $u$-expansions. Let $\sigma$ be an element of $\mathbf{Aut}(L/\FF_q(t,\theta))$ and, 
for $\mu\in M^!_{w,m}\otimes L$, let us denote by $\mu^\sigma\in M^!_{w,m}\otimes L$
the form obtained applying $\sigma$ on every coefficient of the $u$-expansion of $\mu$.

Since $f,\bsb{d}_2$ and $\tau\bsb{d}_2$ are defined over $\FF_q[t,\theta]$,
we get 
$f=\phi^\sigma\bsb{d}_2+\psi^\sigma\tau \bsb{d}_2$, so that 
$$(\phi-\phi^\sigma)\bsb{d}_2+(\psi-\psi^\sigma)\tau\bsb{d}_2=0.$$
Assuming that $\phi^\sigma\neq\phi$ or $\psi^\sigma\neq\psi$ is impossible by Proposition \ref{embedding}. Hence, for all $\sigma$, $\phi^\sigma=\phi$ and $\psi^\sigma=\psi$.
This means that the $u$-expansions of $\phi,\psi$ 
are both defined over $\FF_q(t^{1/q^s},\theta^{1/q^s})$ for some $s\geq 0$. By the fact that
$\mathcal{G}_{\alpha,0,0}=-1+o(1)$ (this follows from the first part of Proposition \ref{interpretationvk}), we get that $\lambda_\alpha\in\FF_q(t^{1/q^s},\theta^{1/q^s})$.

We have proven that $\lambda_\alpha=\widetilde{\pi}^{-\alpha}L(\chi_t,\alpha)(t-\theta)s_\carlitz\in\FF_q(t^{1/q^s},\theta^{1/q^s})$. But we already know that $L(\chi_t,\alpha)\in K_\infty[[t]]$, $s_\carlitz\in K^{\text{sep}}[[t]]$,
and $\widetilde{\pi}\in K_\infty^{\text{sep}}$ (the separable closure of $K_\infty$). Therefore, $$\lambda_\alpha\in\FF_q(t^{1/q^s},\theta^{1/q^s})\cap K_\infty^{\text{sep}}((t))=\FF_q(t,\theta).$$ \CVD

\noindent\emph{Remark.} Proposition \ref{interpretationvk} tells that $\phi\in M_{\alpha+1,1}\otimes\LL$, which is a more precise statement than just saying that
$\phi\in M^!_{\alpha+1,1}\otimes\LL$, following Proposition \ref{embedding}. We also have $\psi=(f-\phi\bsb{d}_2)/\tau\bsb{d}_2$. Since $\bsb{d}_2=1+\ldots$ has $u$-expansion 
in $\FF_q[t,\theta][[u]]$ by the fourth part of Proposition \ref{prarchiv}, the same property holds for
$\tau\bsb{d}_2$ and $\psi$ also belongs to $M_{\alpha+q,1}\otimes\LL$.
We observe the 
additional information that both $\mathcal{G}_{\alpha,0,0}$ and $\tau\mathcal{G}_{\alpha,0,-1}$
are defined over $\FF_q(t^{1/q^s},\theta^{1/q^s})$; this also follows from Proposition \ref{embedding}.



\end{document}

%% file: option_keys
\ifx\optionkeymacros\undefined\else \fi

\catcode`\Œ=\active\defŒ{{\aa}}       
\catcode`\º=\active\defº{\int}        
\catcode`\=\active\def{\c c}        
\catcode`\¶=\active\def¶{\partial}    
\catcode`\Ä=\active\defÄ{\oint}       
\catcode`\Æ=\active\defÆ{\triangle}   
\catcode`\Â=\active\defÂ{\neg}        
\catcode`\µ=\active\defµ{\mu}         
\catcode`\¿=\active\def¿{{\o}}        
\catcode`\¹=\active\def¹{\pi}         
\catcode`\Ï=\active\defÏ{{\oe}}       
\catcode`\§=\active\def§{{\ss}}       
\catcode`\ =\active\def {\dagger}     
\catcode`\Ã=\active\defÃ{\sqrt}       
\catcode`\·=\active\def·{\Sigma}      
\catcode`\Å=\active\defÅ{\approx}     
\catcode`\½=\active\def½{\Omega}      
\catcode`\£=\active\def£{{\it\$}}     
\catcode`\°=\active\def°{\infty}      
\catcode`\¤=\active\def¤{{\S}}        
\catcode`\¦=\active\def¦{{\P}}        
\catcode`\¥=\active\def¥{\bullet}     
\catcode`\»=\active\def»{\leavevmode\raise.585ex\hbox{\b a}}      
\catcode`\¼=\active\def¼{\leavevmode\raise.6ex\hbox{\b o}}        
\catcode`\­=\active\def­{\not=}       
\catcode`\²=\active\def²{\leq}        
\catcode`\³=\active\def³{\geq}        
\catcode`\Ö=\active\defÖ{\div}        
\catcode`\É=\active\defÉ{{\dots}}     
\catcode`\¾=\active\def¾{{\ae}}       
\catcode`\Ç=\active\defÇ{\ll}         
\catcode`\Ò=\active\defÒ{``}          
\catcode`\Á=\active\defÁ{!`}          
\catcode`\¢=\active\def¢{\rlap/c}     
\catcode`\Ô=\active\defÔ{`}           
\catcode`\Õ=\active\defÕ{'}           

\catcode`\=\active\def{{\AA}}       
\catcode`\'=\active\def'{\c C}        
\catcode`\¯=\active\def¯{{\O}}        
\catcode`\¸=\active\def¸{\Pi}         
\catcode`\Î=\active\defÎ{{\OE}}       
\catcode`\®=\active\def®{{\AE}}       
\catcode`\×=\active\def×{\diamond}    
\catcode`\¡=\active\def¡{\accent'27}  
\catcode`\Ó=\active\defÓ{''}          
\catcode`\±=\active\def±{\pm}         
\catcode`\È=\active\defÈ{\gg}         
\catcode`\À=\active\defÀ{?`}          
\catcode`\Ð=\active\defÐ{--}          
\catcode`\Ñ=\active\defÑ{---}         

\catcode`\Š=\active\defŠ{\"a}        
\catcode`\'=\active\def'{\"e}        
\catcode`\•=\active\def•{\"{\i}}     
\catcode`\š=\active\defš{\"o}        
\catcode`\Ÿ=\active\defŸ{\"u}        
\catcode`\Ø=\active\defØ{\"y}        
\catcode`\€=\active\def€{\"A}        
\catcode`\…=\active\def…{\"O}        
\catcode`\†=\active\def†{\"U}        
\catcode`\‡=\active\def‡{\'a}        
\catcode`\Ž=\active\defŽ{\'e}        
\catcode`\'=\active\def'{\'{\i}}     
\catcode`\—=\active\def—{\'o}        
\catcode`\œ=\active\defœ{\'u}        
\catcode`\ƒ=\active\defƒ{\'E}        
\catcode`\ˆ=\active\defˆ{\`a}        
\catcode`\=\active\def{\`e}        
\catcode`\"=\active\def"{\`{\i}}     
\catcode`\˜=\active\def˜{\`o}        
\catcode`\=\active\def{\`u}        
\catcode`\Ë=\active\defË{\`A}        
\catcode`\‹=\active\def‹{\~a}        
\catcode`\–=\active\def–{\~n}        
\catcode`\›=\active\def›{\~o}        
\catcode`\Ì=\active\defÌ{\~A}        
\catcode`\"=\active\def"{\~N}        
\catcode`\Í=\active\defÍ{\~O}        
\catcode`\‰=\active\def‰{\^a}        
\catcode`\=\active\def{\^e}        
\catcode`\"=\active\def"{\^{\i}}     
\catcode`\™=\active\def™{\^o}        
\catcode`\ž=\active\defž{\^u}        

\let\optionkeymacros\null